\documentclass[11pt]{amsart}

\usepackage{amssymb}
\usepackage{tikz}

\newtheorem{theorem}{{Theorem}}[section]
\newtheorem{isom.ext}[theorem]{{Trivial isometric extension}}
\newtheorem{lemma}[theorem]{{Lemma}}
\newtheorem{corollary}[theorem]{{Corollary}}
\newtheorem{remark}[theorem]{{Remark}}
\newtheorem{example}[theorem]{{Example}}

\newcommand{\Z}{{\mathbb{ Z}}}

\begin{document}

\title{Entropy of induced maps of regular curves homeomorphisms}
\author{Aymen Daghar and Issam Naghmouchi}
\address{ Issam Naghmouchi, University of Carthage, faculty
of sciences of bizerte, (UR17ES21), "Dynamical Systems and their Applications"
\\ 7021, Jarzouna, Tunisia}
\email{issam.naghmouchi@fsb.rnu.tn and issam.nagh@gmail.com}

\address{ Aymen Daghar, University of Carthage, faculty
of sciences of bizerte, (UR17ES21), "Dynamical Systems and their Applications"
\\ 7021, Jarzouna, Tunisia}
\email{aymendaghar@gmail.com}
\date{\today}
\maketitle

\begin{abstract}
Let $f:X\to X$ be a self homeomorphism of a continuum $X$, we show that the topological entropy of the induced system $(2^X,2^f)$ is infinite provided that $X\setminus \Omega(f)$ is not empty. If furthermore $X$ is a regular curve then it is shown that $(2^X,2^f)$ has infinite topological entropy if and only if  $X\setminus \Omega(f)$ is not empty. Moreover we prove for the induced system $(C(X),C(f))$ the equivalence between the following properties: (i) zero topological entropy; (ii) there is no Li-Yorke pair and (iii) for any periodic subcontinnum $A$ of $X$ and any connected component $C$ of $X\setminus \Omega(f)$, $C\subset A$ if $A\cap C\neq \emptyset$. In particular, the topological entropy of either $(2^X,2^f)$ or  $(C(X),C(f))$ has only two possible values $0$ or $\infty$. At the end, we give an example of a pointwise periodic rational curve homeomorphism $F:Y\to Y$ with infinite topological entropy induced map $C(F)$.
\end{abstract}


\section{ \bf Introduction}
A \textit{continuum} is a non-empty compact connected metric space. A continuum is a \textit{regular curve} (resp. a \textit{rational curve}) if any point has an $\epsilon$-open neighborhood with finite (resp. at most countable) boundary, for any $\epsilon>0$. Regular curves are known to be $1$-dimensional Peano-continua (see \cite{Kura}). Recall that the class of regular curve is larger than those of arcs, graphs and dendrites.

Let $X$ be a compact metric space. We denote by $2^X$ (resp. $C(X)$) the set of all non-empty compact subsets (resp. compact connected subsets) of $X$. The Hausdorff metric $d_H$ on $2^{X}$ is defined as follows: $d_H (A,B) = \max \Big(\sup_{a\in A} d(a,B),\sup_{b\in B} d(b,A)\Big)$, for each $A, B \in 2^X$ where $d(x,M) = \inf_{y \in M} d(x,y)$ for each $x\in X$ and $M\in 2^X$. With this distance, $(2^X, d_H)$ and $(C(X), d_H)$ are compact metric spaces (for more details, see \cite{conti}). Let $f: X\to X$ be a continuous map of $X$. We denote by $2^{f}:2^{X}\to 2^{X}, A\to f(A)$ and $C(f): C(X) \to C(X),\; C\to f(C)$. $2^{f}$ and $C(f)$ are called induced maps of $f$. Observe that $2^{f}$ (resp. $C(f))$ is a continuous self mapping of $(2^{X},d_{H})$ (resp. $C(X),d_{H})$ (cf. \cite{conti}).

In the past decades, extensive literature has been developed on the dynamical properties of induced maps (see for example \cite{banks}, \cite{Guirao}, \cite{nagar}, \cite{Kwiet}, ect...). In this paper, we will focus on topological entropy of induced maps arising from a wide class of one-dimensional continua homeomorphisms. Lets first give a brief historic review of this topic. Lampart and Raith proved in \cite{Lambard} that the induced map $C(f)$ has zero topological entropy for any interval or circle homeomorphism $f$. For tree homeomorphisms, it can be deduced from  Matviichuk's results in \cite{Treee} that $\omega$-limit sets of the induced map defined on the hyperspace of subcontinua are finite and so there is no Li-Yorke pairs.  Acosta, Illanes and Mendez in \cite{Acosta example} gave an example of a dendrite homeomorphism $f:D\to D$ for which the induced map $C(f)$ has infinite topological entropy (see also \cite{IH} for another example). Recently Hernandez and Mendez in \cite{Hendarez dendrite}, have shown that the topological entropy of the induced maps $C(f)$ and $2^{f}$ of any dendrite homeomorphism $f$ takes only two possible values $0$ of $\infty$. Furthermore, they showed that the existence of a point $x$ lying outside the minimal subdendrite containing both $\omega_f(x)$ and $\alpha_f(x)$ characterizes the case of infinite topological entropy induced map $C(f)$.

In this paper, we will study the topological entropy of induced maps of a given regular curve homeomorphism $f:X\to X$. We show for instance that the topological entropy of both induced maps $2^f$ and $C(f)$ has only two possible values $0$ and $\infty$. Moreover we will give a necessarily and sufficient condition replacing that of Hernandez and Mendez's one (given for dendrites homeomorphisms) which characterizes the case of infinite topological entropy induced map $C(f)$.

\textbf{Plan of the paper}: In section $2$, we give some definitions and preliminary results. In section 3, we will focus on the induced map $2^{f}$. We prove in the first part that for any continuum $X$, $2^{f}$ has infinite topological entropy if $\Omega(f) \subsetneq X$. Then in the second part we deal with regular curves homeomorphisms and we prove that $2^{f}$ has zero topological entropy if and only if $\Omega(f)=X$ (which is known to be equivalent to the equicontinuity of $(X,f)$). In particular, it is shown that topological entropy of $2^f$ has only two possible value $0$ and $\infty$. In section 4, we deal with the induced map $C(f)$ of a regular curve homeomorphism and we give a characterization of the case of infinite topological entropy of the induced map $C(f)$. It is shown also that the topological entropy of $C(f)$ takes only two possible values $\{0,\infty\}$. In section 5, we give an example of a pointwise periodic homeomorphism $F$ of a rational curve with infinite topological entropy induced map $C(F)$.

\section{\bf Definitions and preliminary results}

A dynamical system is a pair $(X, f )$, where $X$ is a compact metric space and $f : X\to X$ is a continuous map. Let $\mathbb{Z},\ \mathbb{Z}_{+},\; \mathbb{Z}_{-}$ and $\mathbb{N}$ denote the sets of integers, non-negative integers, non-positive integers and positive integers, respectively.

 Let $(X, f )$ be a dynamical system.
The $\omega$-limit set of a given point $x\in X$ is defined as follow:
$$\omega_f(x)=\cap_{n\in\mathbb{N}}\overline{\{f^k(x): k\geq n\}}$$

$$=\{y\in X: \exists n_1<n_2<\cdots: \lim_{i \rightarrow +\infty} d(f^{n_i}(x),y)=0\}$$

If $f$ is a homeomorphism then the full (resp. the backward) orbit under $f$ of a given point $x \in X$ is $O_f (x):=\{f^n(x): n\in \Z \}$ (resp. $O^-_f (x):=\{f^n(x): n\in \Z_{-} \}$). The $\alpha$-limit set  of $x$ is $\alpha_f(x):=\omega_{f^{-1}}(x)$.

A point $x$ in $X$ is called
\begin{itemize}
\item \emph{fixed} if $f(x)=x$,
\item \emph{periodic} if $f^n (x)=x$ for some $n\in \mathbb{N}$,
\item \emph{almost periodic}  if for any neighborhood $U$ of $x$ there exists $N\in\mathbb{N}$ such that $\{f^{n+i}(x),\ i=0,1,\cdots,N\} \cap U \neq\emptyset$ for all $n\in\Z_{+}$,
\item  \emph{recurrent} if $x\in \omega_f(x)$,
 \item \emph{regularly recurrent} if for any neighborhood $U$ of $x$ there exists $k\in \mathbb{N}$  such that $f^{kn}(x)\in U$ for all $n\in \mathbb{N}$,
\item \emph{wandering} for if there exists a neighborhood $U$
of $x$ (called wandering neighborhood ) such that $f^{-n}(U)\cap U=\emptyset$, for every $n\in\mathbb{N}$. Otherwise, the point
$x$ is said to be \emph{non-wandering}.
\end{itemize}

We denote respectively by $P(f)$, $AP(f)$, $R(f)$, $RR(f)$  and $\Omega(f)$, the sets of periodic points, almost periodic points, recurrent points, regularly recurrent and the  non-wandering points of $f$. The forward orbit of a given point $x\in X$ is the subset $\{f^n(x): n\in\mathbb{Z}_{+}\}$. A subset $M$ is called \emph{minimal} if it is not empty, closed, $f$-invariant (\textit{i.e.} $f(A)\subset A$) and there is no proper subset of $M$
having these properties.
 Let $A$ be a subset of $X$, we denote by $\Delta_{A}=\{(x,x),\; x\in A\}$. A pair $(x,y)\in X\times X$ is called \textit{proximal} if $\underset{n\to \infty}\liminf\ d(f^n(x),f^n(y))=0$ otherwise it is called \textit{distal}. If $\underset{n\to \infty}\limsup\ d(f^n(x),f^n(y))=0$, then $(x,y)$ is called \textit{asymptotic}. A pair $(x,y)$ is called a \textit{Li-Yorke pair} if it is proximal but not asymptotic. The dynamical system $(X,f)$ is called \textit{distal} (resp. \textit{almost distal}) if it has no proximal pair (resp. no Li-Yorke pair).

 The notion of topological entropy for a continuous self-map of a compact metric space was first introduced by Adler, Konheim and McAndrew in \cite{AKM}. We recall here the equivalent definition formulated by Bowen in \cite{Bowen} and independently by Dinaburg \cite{Dinaburg}
  Let $(X,d)$ be a compact metric space and let $f : X \to X$ be a map. Fix $n \geq 1$ and $\epsilon > 0$. A subset $E$ of $X$ is called
$(n,f,\epsilon)$ separated if for every two different points $x,y\in E$, there exists
$0 \leq j \leq n$ with  $d(f^{j}(x),f^{j}(y))>\epsilon$. Denote by sep$(n, f, \epsilon)$ the maximal possible cardinality of an $(n, f, \epsilon$)-separated set in $X$. Then the \textit{topological entropy} of $f$ is defined by:
$$ h(f)=\displaystyle\lim_{\epsilon \to 0} \limsup_{n\to +\infty} \frac{\ln(sep(n, f, \epsilon))}{n}.$$

Given two dynamical systems $(X,f)$ and $(Y,g)$, by \textit{factor map} we mean a continuous surjective $\pi: X\to Y$ satisfying $\pi \circ f = g \circ \pi$, in this case $(X,f)$ is called an extension of $(Y,g)$ and $(Y,g)$ is called \textit{a factor} of $(X,f)$. Moreover if $\pi$ is an homeomorphism we say that $(X,f)$ and $(Y,g)$ are topologically conjugated. Recall that if $(Y,g)$ is factor of $(X,f)$, then $h(g)\leq h(f)$.

A $\mathbb{Z}$-action on $X$ generated by a self homeomorphism $f$ of $X$ is said to be \textit{equicontinoues} if the family $\{f^{n},\; n\in \mathbb{Z}\}$ is equicontinoues.\\

A family $\{A_{i},\; i\in I\}$ of subset of $X$ is said to be \emph{a null family} if for any $\epsilon>0$ there exist a finite subset $J$ of $I$ such that for any $i\in I\setminus J$, we have $diam(A_{i}) < \epsilon $. It is well known that each pairwise disjoint family of subcontinua of a regular curve is null (see \cite{LFS}).

\begin{lemma}\label{f equi}
 Let $(X,f)$ be an equicontinuous dynamical system.\\ Then $(2^{X},2^{f})$ is equicontinuous, and so is $(C(X),C(f))$.
\end{lemma}

%

\begin{lemma}\label{factor factor}
Let $(X,f)$ be a dynamical system and $(Y,g)$ be a factor of $(X,f)$ thought the factor map $\pi$. Then we have the following assertions:\\
\rm{(i)} $(2^{Y},2^{g})$ is a factor of $(2^{X},2^{f})$ thought the factor map $2^{\pi}$.\\
 \rm{(ii)} If moreover $\pi$ is monotone then $(C(Y),C(g))$ is a factor of $(C(X),C(f))$ thought the factor map $C(\pi)$
\end{lemma}

\begin{proof}
\rm{(i)} The proof of assertion (i) is straightforward.\\
\rm{(ii)} For the proof of (ii) it suffices to show the surjectively of $C(\pi)$. This comes immediately from the fact that $\pi$ is monotone.

\end{proof}

\begin{remark}\rm{
Note that if the factor map $\pi:X\to Y$ is not monotone, then $(C(Y),C(g))$  may be not a factor of $(C(X),C(f))$. For this consider the closed interval $[0,1]$ homeomorphism $f:t\to t^{2}$ and consider the factor system $(Y,g)$ induced by collapsing the closure of the full orbit of $\frac{1}{2}$ to a point. The topological entropy of $(C(Y),C(g))$ is infinite, although the topological entropy of $(C([0,1]),C(f))$ is $0$. Therefore $(C(Y),C(g))$ can not be a factor of $(C([0,1]),C(f))$.}
\end{remark}

Some useful results related to regular curve homeomorphisms from \cite{IAM,Ay,I minimal regular,I2 dynamics RC} are collected in the following Theorem.

\begin{theorem}\label{Thintro1}
Let $f : X\longrightarrow X$ be a homeomorphism of a regular curve $X$ then the following assertions hold:\\
\rm{(1)} $\omega_f(x)$ is a minimal set, for all $x\in X$.\\
\rm{(2)} $\Omega(f) = AP(f)$.\\
\rm{(3)} If $x\in X$ such that $\omega_f(x) \cup \alpha_f (x)$
is infinite then $\omega_f(x) =\alpha_f (x)$.\\
\rm{(4)} The collection of all minimal subsets is closed in $(2^X,d_H)$.\\
\rm{(5)} The map $\omega_f$ is continuous at any point with infinite $\omega$-limit set.\\
\rm{(6)} For any infinite minimal set $M$,  $A(M)=\{x\in X: \omega_f(x)=M\}$ is a closed subset of $X$.\\
\rm{(7)} If $P(f)=\emptyset$, then $\Omega(f)$ is a minimal set and the unique minimal set.\\
\rm{(8)} If $P(f)\neq \emptyset$, then $\Omega(f)=\overline{P(f)}$.\\
\rm{(9)} If $P(f)=\emptyset$, then the proximal relation is a monotone closed equivalence relation and any class meets $\Omega(f)$ in at most two points. Moreover the exists an invariant simple closed curve.
\end{theorem}
In the following, we will give some results further describing the dynamics of regular curves homeomorphisms. Given a regular curve homeomorphism $f:X\to X$, then we have the following:

\begin{lemma}\label{L2,CC asymtotic}
 Every connected component $D$ of $X\setminus \Omega(f)$ satisfies the following two properties:
 \begin{itemize}
  \item [(i)] Every pair of points in $D$ is asymptotic with respect to $f$ and $f^{-1}$. Moreover the $\omega$-limit set (resp. $\alpha$-limit set) of each point in $D$ is the $\omega$-limit set (resp. $\alpha$-limit set) of some point in $\overline{D}\cap \Omega(f)$;
  \item [(ii)] The closure of $D$ intersects $\Omega(f)$ in at most two points.
 \end{itemize}
 \end{lemma}
 \begin{proof}
 (i) There is a proof of (i) in  \cite{{AYYY HLC}}. But we prove it here for the convenience of the reader.  Let $x\in D$ and suppose that $A$ is the intersection of the asymptotic class of $x$ (with respect to $f$) with $D$. As $x$ is a wandering point and $X$ is a regular curve, $A$ is is open in $D$. As each point in $D$ is wandering, then it can be shown easily that $A$ is a closed subset of $D$. It follows by connectedness of $D$ that $A=D$.\\
  For the second part, observe first that $\partial D\subset \Omega(f)$. We distinguish two cases:\\
\textit{Case 1.} $(f^{n}(D))_{n\in \mathbb{Z}}$ is pairwise disjoint: In this case $(f^{n}(D))_{n\in \mathbb{Z}}$ is a null family, thus $\overline{D}\times \overline{D} \subset A(X,f)$. Recall that $\partial D\subset \Omega(f)$. So let $z$ be any point in $\partial D$ then $(z,y)\in A(X,f)\cap A(X,f^{-1})$, $\forall y\in D$ and the result follows.\\
\textit{Case 2.} For some  $p\in \mathbb{N},\; f^{p}(D)=D$: In this case, each point in $D$ is proximal (and so asymptotic) with respect to $f$ (resp. $f^{-1}$) to a periodic point $t$ (resp. $t^{'}$). Hence $ \omega_{f}(y)=O_{f}(t)$ and $\alpha_f(y)=O_f(t^{'})$ for all $y\in D$. Finally $t,t^{'}\in \partial D$ since $f^{p}(\overline{D})=\overline{D}$.

 (ii) If $P(f)=\emptyset$, the proof of property (ii) comes immediately from property (i) and assertion (9) in Theorem \ref{Thintro1}.
  Assume now that $P(f)\neq \emptyset$. If the family $\{f^{n}(D), n\in \mathbb{Z}\}$ is pairwise disjoint then it a null family and so is  $\{f^{n}(\overline{D}), n\in \mathbb{Z}\}$. In this case, $\Omega(f) \cap \overline{D}$ should be reduced to a single point since $(\Omega(f),f_{\mid \Omega(f)})$ is a distal system (see [Theorem 5.2 in \cite{IAM}]). Now if the family $\{f^{n}(D), n\in \mathbb{Z}\}$ is not pairwise disjoint then $f^{k}(D)=D$  (and so $f^{k}(\overline{D})=\overline{D}$) for some $k\in \mathbb{N}$.
   Recall that by property (i) any pair of points in $D$ is asymptotic. Since $f^{k}(D)=D$ , then for some $a,b\in Fix(f^{k})\cap \overline{D},\; \omega_{f^{k}}(x)=\{a\}$, (resp. $\alpha_{f^{k}}(x)=\{b\})\; \forall x\in D$. Suppose that there exists $c\in \overline{D}\cap \Omega(f) \setminus\{a,b\}$. Then obviously $\{a,b\}\cap \omega_{f^{k}}(c)=\emptyset$ (see assertion (1) in Theorem \ref{Thintro1}). Let $U$ be an open neighbourhood of $\omega_{f}(c)$ in $\overline{D}$ with finite boundary such that $\{a,b\}\cap \overline{U}=\emptyset$. By [Theorem 5 in \cite{accecipilite}] the subset $\{c\}\cup D$ is arcwise connected and hence it is locally arcwise connected (see [Corollary 5.5 \cite{CPA iml LCPA}]). Therefore we may find for each $n\geq 0$ an arc $I_{n}$ joining $c$ and a point $c_{n}\in D\cap U$ such that $I_{n}\cap \Omega(f)=\{c\}$ and $\lim_{n \to +\infty} diam(I_{n})=0$. Recall that for any $n\geq 0, \omega_{f^{k}}(c_{n})=\{a\}$ and therefore for each $n\geq 0$, $f^{ks_{n}}(I_{n}) \cap \partial U \neq \emptyset$ for some $s_{n}\geq 0$. It follows that $c\in \alpha_{f^{k}}(x)$ for some point $x\in \partial U\cap D$. This is a contradiction with the fact that the $\alpha$-limit set of any point in $D$ is $\{b\}$.
\end{proof}

\begin{lemma}\label{2 a 2 disjoin}
For each connected component $D$ of $X\setminus \Omega(f)$ for which the closure intersects $\Omega(f)$ in a non-periodic point, the family $\{f^{n}(\overline{D}), \; n\in \mathbb{Z})\}$ is pairwise disjoint. In particular, it is a null family.\\
\end{lemma}

\begin{proof}
 Let $a\in \overline{D}\cap (\Omega(f) \setminus P(f))$. Recall that $\omega_{f}(a)$ is an infinite minimal set (see assertion (2) of Theorem \ref{Thintro1}). We will show that, $\omega_{f}(x)=\omega_{f}(a)$, for any $x\in D$. Indeed, pick any sequence $(a_{n})_{n\geq 0}$ in $D$ with limit $a$. Recall that $(\omega_{f}(a_{n}))_{n\geq 0}$ is constant (see Lemma \ref{L2,CC asymtotic}). It follows from the continuity of the map $\omega_{f}$ at $a$ (see assertion (5) in Theorem \ref{Thintro1}) that $\omega_{f}(a_{n})=\omega_{f}(a)$, for each $n\geq 0$.  Again by Lemma \ref{L2,CC asymtotic}, $\omega_{f}(x)=\omega_{f}(a),\; \forall x\in D$.

 Suppose that $\{f^{n}(\overline{D}), \; n\in \mathbb{Z}\}$ is not pairwise disjoint, then $f^{k}(D)=D$ for some $k\in \mathbb{N}$. Thus $\{z,f^{k}(z)\} \subset D$ for some $z\in D$. By Lemma \ref{L2,CC asymtotic}, the pair $(z,f^{k}(z))$ is asymptotic and so $\omega_{f}(z)$ contains a periodic point. It turns out that $\omega_{f}(a)=\omega_{f}(z)$ contains a periodic point, which is a contradiction. In conclusion, $\{f^{n}(\overline{D}), \; n\in \mathbb{Z})\}$ is pairwise disjoint.

\end{proof}

\begin{lemma}\label{single point}
Suppose that $P(f)\neq \emptyset$ and $D$ is a connected component of $X\setminus \Omega(f)$ for which the closure intersects $\Omega(f)$ in a non-periodic point. Then $\overline{D}\cap \Omega(f)$ is reduced to a single point.
\end{lemma}

\begin{proof}
 By Lemma \ref{2 a 2 disjoin}, any pairs in $\overline{D}$ is an asymptotic pair. Recall also that any pair of distinct points in $\Omega(f)$ is distal. Hence $\overline{D}\cap \Omega(f)$ is reduced to a single point.
\end{proof}

The following lemma follows immediately from the definition of the distance $d_{H}$.
\begin{lemma}\label{renuion et d(H)}
  Given a compact metric space $X,\; A_{1},\dots A_{n},\; B_{1},\dots B_{n}$ a finite family of closed subsets of $X$ and $\epsilon>0$ such that for any $1\leq i\leq n$ we have $d_{H}(A_{i},B_{i})\leq \epsilon$. Then $d_{H}(\displaystyle\bigcup_{1\leq i\leq n}A_{i},\displaystyle\bigcup_{1\leq i\leq n}B_{i})\leq \epsilon$.
\end{lemma}

\section{\bf On the induced map $2^{f}$ of regular curves homeomorphism $f$.}

\begin{theorem}\label{non deno}
 Let $X$ be a compact metric space and $f:X\to X$ be a homeomorphism such that $X\setminus \Omega(f)$ is uncountable, then $h_{top}(2^{f})=+\infty$.
\end{theorem}
For the proof of this theorem, we will use the following lemma.

\begin{lemma}\label{nombre fini}
  Let $X$ be a compact metric space and $f:X\to X$ be a homeomorphism. If for some $n\in \mathbb{N}$, there exists a sequence of points $x_{1},\dots,x_{n}$ satisfying the following conditions:\\
  \rm{(i)} $O_f(x_{i}) \cap O_f(x_{j})=\emptyset$, if $i\neq j$;\\
  \rm{(ii)}$\{x_{1},\dots,x_{n}\} \cap \big( \displaystyle\bigcup_{1\leq i\leq n}\alpha_{f}(x_{i})\cup \omega_{f}(x_{i})\big)=\emptyset$,\\

  then $h_{top}(2^{f})\geq n\ln(2)$.

\end{lemma}

\begin{proof}

Let $Y=\displaystyle\bigcup_{1\leq i\leq n}  \overline{O_{f}(x_{i})}$ and let $\delta=\inf_{1\leq i\leq n}\{d(x_{i},Y\setminus\{x_{i}\})\}$. Observe first that $\delta>0$. Let $0<\epsilon < \delta$ and $k \in \mathbb{N}$.
For each $\sigma=(\sigma_{0},\dots, \sigma_{k-1}) \in \mathcal{P}(\{1,\dots,n\})^{k}$, we define the closed subset $A_{\sigma}$ as follows:
 $$A_{\sigma}=\big(\displaystyle\bigcup_{1\leq i\leq  n} \alpha_{f}(x_{i}) \cup \{f^{-j}(x_{i}), j\geq k\}\big) \cup \{f^{-j}(x_{i}),\; i\in \sigma_{j},\; 0\leq j\leq k-1\}.$$
 We will show that $\{A_{\sigma}, \; \sigma \in \mathcal{P}(\{1,\dots,n\})^{k}\}$ is a $(k,2^f,\epsilon)$-separated set for $(2^{X},2^{f})$. Indeed if $\sigma=(\sigma_{0},\dots, \sigma_{k-1})$ and $\sigma^{\prime}=(\sigma^{\prime}_{0},\dots, \sigma^{\prime}_{k-1})$ are two distinct elements of $\mathcal{P}(\{1,\dots,n\})^{k}$, then $\sigma_{j}\neq \sigma_{j}^{\prime}$ for some $0\leq j\leq k-1$. We may assume without loss of generality that $i\in \sigma_{j} \setminus \sigma_{j}^{\prime}$ for some $1\leq i\leq n$. Thus $x_{i}\in f^{j}(A_{\sigma})\setminus f^{j}(A_{\sigma^{\prime}})$. It follows that $d_{H}(f^{j}(A_{\sigma}),f^{j}(A_{\sigma^{\prime}}))>\epsilon $. In result $\{A_{\sigma}, \; \sigma \in \mathcal{P}(\{1,\dots,n\})^{k}\}$ is a $(k,\epsilon)$-separated set for $(2^{X},2^{f})$. Consequently, $\limsup_{k\to +\infty}\frac{\ln(sep(n, f, \epsilon))}{k}\geq n\ln(2)$ and so $h_{top}(2^{f})\geq n\ln(2)$.
\end{proof}

\textit{Proof of Theorem \ref{non deno}}.
  As $X\setminus \Omega(f)$ is an invariant uncountable subset, then for each $n\in \mathbb{N}$, we may find $n$ distinct points $x_{1},\dots,x_{n}$ satisfying conditions $(i)$ and $(ii)$ of Lemma \ref{nombre fini}. Therefore $h_{top}(2^{f})=+\infty$.

 \qed

\begin{corollary}\label{Omega=X infini 2f}
Let $X$ be a continuum and $f:X\to X$ be a homeomorphism such that $X\setminus \Omega(f) \neq \emptyset$, then $h_{top}(2^{f})=+\infty$.
\end{corollary}
\begin{proof}
By boundary bumping Theorem (see \cite{conti}, 5.6, page 74), any connected component of $X\setminus \Omega(f)$ is not degenerate thus uncountable. It turns out that $X\setminus \Omega(f)$ is uncountable. By Theorem \ref{non deno}, $h_{top}(2^{f})=+\infty$.
\end{proof}

\begin{corollary}\rm{
Let $f:X\to X$ be a regular curve homeomorphism. If $X\setminus \Omega(f) \neq \emptyset$, then $h_{top}(2^{f})=+\infty$, otherwise $h(2^{f})=0$. In particular $h_{top}(2^{f})\in \{0,+\infty\}$.}
\end{corollary}
\begin{proof}

Assuming that $\Omega(f)=X$. According to Theorem 5.2 in \cite{IAM}, $(X,f)$ is equicontinuous if $P(f)\neq \emptyset$. If $P(f)=\emptyset$, then by Corollary 4.6 in \cite{Ay} $X$ is a simple closed curve and $f$ is an irrational rotation (up to a conjugation). Thus $(X,f)$ is equicontinoues. By Lemma \ref{f equi}, $(2^{X},2^{f})$ is equicontinoues and so $h(2^{f})=0$.

Assuming now that $X\setminus \Omega(f) \neq \emptyset$, then $h_{top}(2^{f})=+\infty$ by Corollary \ref{Omega=X infini 2f}.
\end{proof}
\begin{remark} \rm{
Note that for a regular curve homeomorphism $f:X\to X$, topological entropy of the induced map $2^{f}$ is infinite if and only if $X\setminus \Omega(f)\neq \emptyset$. This is no more longer true for rational curves homeomorphisms, in fact we will give in \ref{example rational} an example of a pointwise periodic rational curve homeomorphism $F:Y\to Y$ such that $h_{top}(C(F))=h_{top}(2^{F})=+\infty$.                              }
\end{remark}

\section{\bf Topological entropy of the induced map $C(f)$ for a regular curve homeomorphism $f$}

In this section, we give sufficient and necessarily conditions of a regular curve homeomorphism $f$ to have induced map $C(f)$ with zero topological entropy. In particular, we show that for any regular curve homeomorphism $f,\; h_{top}(C(f))\in \{0,+\infty\}$. In the sequel, we will denote by $f:X\to X$ a regular curve homeomorphism.

\begin{theorem}\label{Main 1 eq} The following properties are equivalent: \\
(i) The topological entropy of $C(f)$ is zero;\\
(ii) $(C(X),C(f))$ has no Li-Yorke pair;\\
(iii) For any periodic subcontinuum $A$ of $X$ and any connected component $C$ of $X\setminus \Omega(f)$,  $C\subset A$ whenever $A\cap C\neq \emptyset$;\\
(iv) For any periodic subcontinua $A$ and $B$ of $X$ such that $A\subset B$, any connected component of $B\setminus A$ has finite orbit.\\
If moreover $X$ is different from a simple closed curve then the above properties are all equivalent to the following:\\
(v) Any subcontinuum of $X$ is asymptotic to a regularly recurrent one.
\end{theorem}

\begin{remark}
\rm{
The condition given in \cite{Hendarez dendrite} characterizing the infinitude of the topological entropy of the induced map $C(f)$ of a dendrite homeomorphism $f$ could not be stated in the case of regular curve, since regular curves are in general not uniquely arcwise connected. Needless to say, the opposite of property (iv) given in Theorem \ref{Main 1 eq} is equivalent to the one given in \cite{Hendarez dendrite} in the case of dendrite. Indeed, suppose that there exists a point $x\in X$ which does not belong to the minimal subdendrite $D_{min}$ containing $\omega_{f}(x)$ and $\alpha_{f}(x)$. Denote by $B=D_{min}\cup \cup_{n\in\mathbb{Z}} f^n(I)$ where $I$ is the arc joining $x$ to a point $y\in D_{min}$ such that $I\cap D=\{y\}$. We let $A=D_{min}$ then we can check easily the opposite of property (iv). Now suppose that we have the opposite of property (iv). Let $C$ be a connected of $B\setminus A$ with infinite orbit. Thus $\lim_{\mid n\mid\to\infty}diam(f^n(\overline{C}))=0$ and so the $\omega$-limit set (as well as the $\alpha$-limit set) of any point in $C$ is a subset of $\cup_{N\in\mathbb{Z}} f^n(A)$.  Take a point $x\in C$ and denote by $p$ a common period of $B$ and $A$. Let $D$ be the convex hull of the orbit of $A$. Then $D$ consists of the union of the orbit of $A$ and the union of the  orbit of a periodic arc joining a periodic point of $A$ to a periodic point of some $f^s(A)$ with $f^s(A)\neq A$. Clearly $x\notin D$ and the convex hull of its $\omega$-limit and the $\alpha$-limit set is included into $D$.
}
\end{remark}

For the proof of Theorem of \ref{Main 1 eq}, we need the following lemmas.

\

\begin{lemma}\label{L1,premier equiv}
If there exists a periodic subcontinuum $A$ intersecting but not containing a connected component $C$ of $X\setminus \Omega(f)$, then there exist a periodic subcontinuum $B$ of $X$ such that $A \subset B$ and some connected component of $B\setminus A$ with infinite orbit.
\end{lemma}

\begin{proof}
Without loss of generality, we may assume that $f(A)=A$. Observe that $T=A\cap C$ is a proper closed subset of $C$. By connectedness of $C, T$ is not open in $C$, so there exists $x\in T$ such that $x\notin int_{C}(T)$. Recall that any point of $C$ is wandering, then let $O$ be an open connected neighbourhood of $x$ in $C$ (in particular in $X$) such that the family $\{f^{n}(O),\; n\in \mathbb{Z}\}$ is pairwise disjoint. As $X$ is a regular curve, $\{f^{n}(O),\; n\in \mathbb{Z}\}$ is a null family. Pick an arc $I\subset O$ joining a point $y\in O\setminus A$ and a point $z\in A$ and such that $I\cap A=\{z\}$. Note that for any $n\in \mathbb{Z}$, the arc $f^{n}(I)$ also intersects $A$ in only one point which is $f^{n}(z)$. Moreover $\{f^{n}(I),\; n\in \mathbb{Z}\}$ is a null family. Let $B=A\cup (\cup_{n\in \mathbb{Z}} f^n(I))$. Thus $B$ is a fixed point of $(C(X),C(f))$ and $I$ is a connected component of $B\setminus A$ having infinite orbit.
\end{proof}

\begin{lemma}\label{L3,cc arc ou cercle}
Assuming that $X$ is not a simple closed curve and for any periodic subcontinuum $A$ of $X$ and any connected component $C$ of $X\setminus \Omega(f)$, $C\subset A$ whenever $A\cap C\neq \emptyset$. Then the following hold for any connected component $D$ of $X\setminus \Omega(f)$:
\begin{itemize}
   \item [(i)]  $\overline{D}\cap \Omega(f) \subset P(f)$;
  \item [(ii)] $\overline{D}$ is a simple closed curve if $\overline{D}\cap \Omega(f)$ is reduced to a single point, otherwise it is an arc.
  \item[(iii)] $f^{k}(D)=D$, for some $k\in \mathbb{N}$
\end{itemize}
\end{lemma}
\begin{proof}

 (i) First we will prove that $P(f)\neq \emptyset$. Suppose that $P(f)=\emptyset$. By assertion (9) in Theorem \ref{Thintro1}, we may find $S \subsetneq X$ an invariant simple closed curve. Let $x\in X\setminus S$ and let $I$ be an arc joining $x$ to a point $y$ in $S$ such that $I\cap S=\{y\}$. Since $\Omega(f) \subset S$, then from Lemma \ref{2 a 2 disjoin} that $\{f^{n}(I), n\in \mathbb{Z}\}$ is a null family. Let $x_{0}\in I\setminus \{x,y\}$ and let $I_{0}$ be the subarc of $I$ joining $x_{0}$ and $y$. Set $A=S\cup \big(\displaystyle\bigcup_{n\in \mathbb{Z}}f^{n}(I_{0})\big)$. Clearly $A$ is a fixed subcontinuum of $X$. Let $C$ be the connected component of $X\setminus \Omega(f)$ containing $x$. Then $C\cap A\neq \emptyset$ and $C \nsubseteq A$, a contradiction with the assumption of the lemma. Therefore  $P(f)$ should be not empty.

 Let $D$ be a connected component of $X\setminus \Omega(f)$. Suppose that $\overline{D}\cap \Omega(f)  \nsubseteq P(f)$. Let $a\in \overline{D}\cap \Omega(f) \setminus P(f)$. By Lemma \ref{single point}, $\{a\}=\overline{D}\cap \Omega(f)=\partial D$. Let $Y=X\setminus (\displaystyle\bigcup_{n\in \mathbb{Z}}f^{n}(D))$, obviously $Y$ is a closed strongly invariant subset of $X$. We will show moreover that $Y$ is arcwise connected. Indeed, let $u$ and $v$ be two distinct points in  $Y$. Recall that $\partial (f^{n}(D))=\{f^{n}(a)\}$, for each $n\in \mathbb{Z}$. Thus any arc in $X$ joining $u$ and $v$ is entirely included into $Y$. Consequently $Y$ is arcwise connected. Let $L$ be an arc in $\overline{D}$ joining $a$ to some point $b\in D$ such that $L \subsetneq \overline{D}$. Let $A=Y \bigcup (\displaystyle\bigcup_{n\in \mathbb{Z}} f^{n}(L))$. It is easy to see that $A$ arcwise connected and strongly invariant. Furthermore $A$ is closed by Lemma \ref{2 a 2 disjoin}. Therefore $A$ is a periodic subcontinuum of $X$, moreover $A\cap D \neq \emptyset$ and $D \nsubseteq A$, a contradiction with the assumption of the lemma. We conclude then that $\overline{D}\cap \Omega(f)  \subset P(f)$.\\

(ii)  By assertion (i), $\overline{D}\cap \Omega(f) \subset P(f)$. Moreover by Lemma \ref{L2,CC asymtotic}, $\overline{D}\cap \Omega(f)$ contains at most two points. So let $N\geq 0$ and $g=f^{N}$ such that $\overline{D}\cap \Omega(g) \subset Fix(g)$.\\
\textit{Claim 1}; For some $k\geq 0, g^{k}(D)=D$. If not then the family $\{g^{k}(D),k\in \mathbb{Z}\}$ is pairwise disjoint and so it is null. We let then $A=\displaystyle\bigcup_{n\in \mathbb{Z}}g^{n}(I)$, where $I$ is an arc in $D$ joining a point $x\in \overline{D}\cap \Omega(g)$ to a point $y\in D$ such that $I \subsetneq D$. Clearly $A$ is a periodic subcontinuum of $X$, moreover $A\cap D\neq \emptyset$ and $D\nsubseteq A$, a contradiction with the assumption of the lemma. Thus $g^{k}(D)=D$ for some $k\in \mathbb{N}$

Let $x\in D$ and let $I$ be an arc in $D$ joining $x$ to $g^{k}(x)$. Clearly $I$ contains only finitely many points of the orbit of $x$ under $g^{k}$ (otherwise $I\cap \Omega(g)\neq \emptyset$). So we may assume that $I$ meets the orbit of $x$ only at $x$ and $g^{k}(x)$ (if not consider the subarc $J$ of $I$ joining $x$ and $g^{ks}(x)$ for some $s\in \mathbb{Z^{*}}$ such that $J\cap O_{g^{k}}(x)=\{x,g^{ks}(x)\}$) and consider $g^{s}$ instead of $g$.

\textit{Claim 2} $\{g^{ks}(I), \; s\in \mathbb{Z}\}$ is a null family. We have $I\subset X\setminus \Omega(g)$, thus we can find a finite open cover $\{O_{1},\dots O_{n}\} $ of $I$ in $X$ such that for each $1\leq i\leq n,\; \{g^{ks}(O_{i}), \; s\in \mathbb{Z}\}$ is a null family. Since $I$ is an arc, we may assume that $O_{i}\cap O_{i+1}\neq \emptyset$, for each $1\leq i\leq n-1$. For any $s\in \mathbb{Z}$, we have $g^{ks}(I) \subset g^{ks}(O_{1})\cup\dots \cup f^{ks}(O_{n})$ and thus $diam(g^{ks}(I))\leq diam(f^{ks}(O_{1}))+\dots diam(f^{ks}(O_{n}))$. It follows that $\{g^{ks}(I), \; s\in \mathbb{Z}\}$ is a null family.

Let $J=\overline{\displaystyle\bigcup_{n\in \mathbb{Z}}g^{kn}(I)}$. By claim 2, we have $J=\displaystyle\bigcup_{n\in \mathbb{Z}}g^{kn}(I)\cup \omega_{g^{k}}(x)\cup \alpha_{g^{k}}(x)$, where $\omega_{g^{k}}(x)$ (resp. $\alpha_{g^{k}}(x)$) is reduced to a single fixed point under $g^{k}$ contained in $\overline{D}\cap \Omega(g)$.
 Clearly $J$ is a periodic subcontinuum of $X$ meeting $D$. By assumption of the lemma $D\subset J$, thus $\overline{D}=J$ and $D=\displaystyle\bigcup_{n\in \mathbb{Z}}g^{kn}(I)$.

 \textit{Claim 3}. Any point of $\overline{D}$ is of order less or equal to $2$.
 First observe that for any $s<s^{\prime}\in \mathbb{Z}$

\begin{eqnarray}
g^{ks}(I)\cap g^{ks^{\prime}}(I)= \left\lbrace\begin{array}{ccc}
           \emptyset    &\mbox{if}\; s^{\prime}> s+1 \\
          g^{k(s+1)}(x) &\mbox{if}\; s^{\prime}= s+1
             \end{array}\right.
 \end{eqnarray}

Indeed if not then we can find $a\in I \setminus \{x\}$ and $q\in \mathbb{Z}$ such that $|q|\geq 1$ and $g^{kq}(a)\in I$. Observe that $\{a,g^{kq}(a)\} \subset I\setminus \{x,g^{k}(x)\}$. By the same arguments as above, we prove that $D=\displaystyle\bigcup_{m\in \mathbb{Z}}g^{kqm}(I_{a})$, where $I_{a}$ is the subarc of $I$ joining $a$ and $g^{ks}(a)$. Recall that $O_{g^{k}}(x) \subset D$, therefore $O_{g^{k}}(x)$ should meet $I_{a}$. Since $I_{a}\subset I\setminus \{x,g^{k}(x)\}$, we get a contradiction with the fact that $I$ contains only two points from the orbit of $x$ under $g^{k}$. This implies (1). Then any point in $D$ has an open arc as a neighbourhood in $\overline{D}$ and thus any point in $D$ is of order $2$ in $\overline{D}$. It remains to show that any point in $\partial D$ is of order equal or less to $2$. Take $z\in \partial D$ and assume that $z$ is of order greater or equal to $3$. By the $n$-Beinsatz Theorem (see \cite{Kura}, page 277), we may find a $3$-star $S$ centred at $z$ and included in $\overline{D}$. By Lemma \ref{L2,CC asymtotic}, $\partial U$ is finite. Then take a $3$-star $T\subset S$ small enough such that $T\setminus \{z\} \subsetneq D$. Recall that $D$ is arcwise connected, so let $c\in D\setminus T$ and let $L_{1},L_{2}$ and $L_{3}$ be three arcs in $D$ each of which join $c$ and an end point of $T$. This will form a $3$-star $T^{\prime}\subset L_{1}\cup L_{2}\cup L_{3} \subset D$, this is a contradiction since all points of $D$ are of order $2$. In conclusion, every point of $\partial D$ is of order less or equal to $2$. By [Theorem 9.5 in \cite{conti}], $\overline{D}$ is either an arc or a simple closed curve.

 Recall that $\overline{D}\setminus (\overline{D}\cap \Omega(g))=D$ is connected. Therefore if $\overline{D}\cap \Omega(g)$ contains exactly two points, then $\overline{D}$ should be an arc ( otherwise $\overline{D}$ is a simple closed curve and $\overline{D}\setminus \{u,v\}$ is not connected for any $u,v\in \overline{D}$. In particular, $D=\overline{D}\setminus (\overline{D}\cap \Omega(g))$ is not connected, a contradiction).

 Now if $\overline{D}\cap \Omega(g)$ is reduced to a single point. Then $\overline{D}$ should be a simple closed curve (otherwise $\overline{D}$ is an arc and one of its end points belongs to $D$ and thus it is of order $2$ in $\overline{D}$, a contradiction).

 (iii) We will assume that $f^{k}(D)\neq D$ for any $k\in \mathbb{N}$. In this case $\{f^{k}(\overline{D}),\; k\in \mathbb{Z}\}$ is a null family. By (ii) $\partial D \subset P(f)$ and contains at most two point and since $\overline{D}$ is arcwise connected, let $I \subsetneq D \cup \{a\}$ bean arc joining a point $a\in \partial D$ and a point $x\in D$. We may easily check that $A=\displaystyle\bigcup_{n\in \mathbb{Z}}f^{pn}(I)$, is a periodic subcontinuum of $X$, where $p$ is the period of $a$. Furthermore $A$ meets $D$ but does not contain it, this contradicts the assumptions of the Lemma.

\end{proof}

\begin{lemma}\label{attracteur}
 Assuming that $X$ is not a simple closed curve and for any periodic subcontinuum $A$ of $X$ and any connected component $C$ of $X\setminus \Omega(f)$, $C\subset A$ whenever $A\cap C\neq \emptyset$.  Let $D$ be a connected component of $X\setminus \Omega(f)$ of period $k$. Let $L$ be a subcontinuum of $X$ having a non-empty intersection with $D$. Then the sequence $(f^{kn}(L\cap \overline{D}))_{n\geq 0}$ converges to either $\overline{D}$ or some subset of $\partial D$.
   \end{lemma}
\begin{proof}
By Lemma \ref{L3,cc arc ou cercle}, $\overline{D}$ is either a periodic arc or a periodic simple closed curve. If $\overline{D}$ is an arc with endpoints $a$ and $b$, then $a,b\in Fix(f^{k})$. We can assume that $a$ is an attracting fixed point of $(\overline{D},f_{\mid \overline{D}}^{k})$ and $b$ is a repelling point. If $b$ belongs to a non degenerated connected component of $L\cap \overline{D}$, then $(f^{kn}(L\cap \overline{D}))_{n\geq 0}$ converges to $\overline{D}$. Otherwise the sequence $(f^{kn}(L\cap \overline{D}))_{n\geq 0}$ converges to either $\{a\}$ or $\{a,b\}$.
Now If $\overline{D}$ is a simple closed curve, with $\partial D=\{a\}\subset Fix(f^{k})$. Then $\overline{D}\cap L$ is a subcontinuum of $\overline{D}$ containing $a$. Choose an orientation of $\overline{D}$, so that for any $x\in \overline{D}\setminus \{a\},\; f^{k}([a,x]) \supsetneq [a,x]$ (where $[a,x]$ is the sub arc of $\overline{D}$ joining $a$ and $x$ with respect to this orientation.) Now if $L\cap \overline{D} \supseteq [a,x]$ for some $x\in \overline{D}\setminus\{a\}$, then The sequence $(f^{kn}(L\cap \overline{D}))_{n\geq 0}$ converges to $\overline{D}$, otherwise it converges to $\{a\}$.
\end{proof}

Let $L$ be a subcontinuum of $X$ and $D$ be a connected component of $X\setminus \Omega(f)$. We will say that $L$ satisfies condition $\mathcal{P}(D)$ if $L\cap D \neq \emptyset$ and $\displaystyle\lim_{n\to +\infty}d_{H}(f^{kn}(L\cap \overline{D}),\overline{D}))=0$, where $k$ is the period of $\overline{D}$.

\begin{remark}\label{att} \rm{
   By Lemma \ref{attracteur}, if $L$ does not satisfies condition $\mathcal{P}(D)$, for some connected component $D$ of $X\setminus \Omega(f)$, then either $L\cap D=\emptyset$ or $D\cap L \neq \emptyset$ and $\displaystyle\lim_{n\to +\infty}d_{H}(f^{kn}(L\cap \overline{D}),B)=0$, for some $B\subset \partial D$. This properties will be used in the proof of the lemma below.}
\end{remark}

\begin{lemma}\label{L44}
Assuming that $X$ is not a simple closed curve and for any periodic subcontinuum $A$ of $X$ and any connected component $C$ of $X\setminus \Omega(f)$, $C\subset A$ whenever $A\cap C\neq \emptyset$. Then any subcontinuum of $X$ is asymptotic to a regularly recurrent one. In particular, $(C(X),C(f))$ has no Li-Yorke pairs.
\end{lemma}

\begin{proof}
Let $\{C_n: n\in\zeta\}$ (where $\zeta$ is a subset of $\mathbb{N}$) be the family of connected components of $X\setminus \Omega(f)$ so that the elements of each orbit are placed successively in this indexation. By this way of indexation, if $\zeta=\mathbb{N}$ and $\epsilon>0$ then the diameter of each component as well as the diameter of each element from its orbit is less than $\epsilon$ eventually . Let $A\in C(X)$. We will distinguish several cases depending on the trace of $A$ on $\Omega(f)$. Denote by $A_{\Omega}=A\cap \Omega(f)$.

\textit{Case 1.}  $A_{\Omega}= \emptyset$. Then $A\subset D$ for some connected component of $X\setminus \Omega(f)$. By Lemma \ref{L3,cc arc ou cercle}, $A$ is an arc asymptotic to some periodic single point $\{a\}$ in $\partial D$.

\textit{Case 2.} $A_{\Omega} \neq\emptyset$. For each $n\in\zeta$, let $A_n=\overline{C_n}$ if $A$ satisfies condition $\mathcal{P}(C_n)$ otherwise $A_n=\emptyset$. Denote by $B=A_{\Omega}\cup (\cup_{n\in\zeta}A_n)$.

Claim 1. $B$ is a regularly recurrent point of $(C(X),C(f))$. First, we will verify that $B$ is an element of $C(X)$. Indeed, $B$ is the union of $A_{\Omega}$ with a null family of subcontinua meeting the closed set $A_{\Omega}$, so $B$ is clearly closed. Let $u,v\in A_{\Omega}\subset A$ and let $I$ be an arc of $A$ joining $u$ and $v$. If $I \subset A_{\Omega}$ we are done. Otherwise by Lemma \ref{L3,cc arc ou cercle}, any connected component $C$ of  $I\setminus A_{\Omega}$ is also a connected component of $X\setminus \Omega(f)$, thus $A$ satisfies $\mathcal{P}(C)$ and so $\overline{C}\subset B$. In conclusion, $I\subset B$.

Recall that $\Omega(f)=RR(f)$ and $(\Omega(f),f_{\mid \Omega(f)})$ is an equicontinuous system. It follows that $A_{\Omega}$ is a regularly recurrent point of $(2^X,2^f)$. If $\zeta$ is finite then clearly $B$ is a regularly recurrent point of $(C(X),C(f))$. Assume now that $\zeta=\mathbb{N}$. Given $\epsilon>0$. Let $N\in\mathbb{N}$ be such that for any integer $n> N$, the diameter of each component $C_n$ as well as the diameter of each element from its orbit is less than $\epsilon$. Recall that for some $k\in\mathbb{N}$, $d_H(f^{kl}(A_{\Omega}),A_{\Omega})<\epsilon$ for all $l\in\mathbb{N}$. By Lemma \ref{L3,cc arc ou cercle}, the closure of each component of $X\setminus\Omega(f)$ is periodic. If $m$ is a common period of each $\overline{C_i}, i=1,\dots,N$ then we have

$$d_H(f^{lkm}(B),B)<\epsilon, \forall l\in\mathbb{N}.$$

Claim 2. The pair $(A,B)$ is asymptotic. We will consider only the case $\zeta=\mathbb{N}$, the proof in the case $\zeta$ is finite can be done similarly. Let $\epsilon>0$ and let $N\in\mathbb{N}$ be such that for any integer $n> N$, the diameter of each component $C_n$ as well as the diameter of each element from its orbit is less than $\epsilon$.
By Lemma \ref{attracteur}, there exists $l\in\mathbb{N}$ so that for any $i=1,\dots, N$ such that $A\cap \overline{C_i}\neq\emptyset$, we have for any integer $s>l$

$
\left\{
  \begin{array}{ll}
    d_H(f^{s}(A\cap \overline{C_i}), f^s(\overline{C_i}))<\epsilon, & \hbox{if $A$ satisfies $\mathcal{P}(C_i)$,} \\
    d_H(f^{s}(A\cap \overline{C_i}), f^s(A\cap \partial C_i))<\epsilon, & \hbox{otherwise.}
  \end{array}
\right.$

It turns out by Lemma \ref{renuion et d(H)} that $d_H(f^{s}(A),f^{s}(B))<\epsilon$ for any integer $s>l$.

\end{proof}

Now we will study separately the case of $X$ being a simple closed curve.

\begin{lemma}\label{cercle}
If $X$ is a simple closed curve, then $(C(X),C(f))$ has no Li-Yorke pair. In particular, $(C(X),C(f))$ has zero topological entropy.
\end{lemma}

\begin{proof}
\textit{Case.1} $f$ preserve the orientation: We consider on $X$ the counterclockwise sense, let $I=[a,b]$ and $J=[c,d]$ be two arcs in $X$ such that $(I,J)$ is a proximal pair for $(C(X),C(f))$. Then $(a,c)$ and $(b,d)$ should be proximal pairs for $(X,f)$, thus they are asymptotic pairs for $(X,f)$ (see Lemma 43 in \cite{blanch}). It turns out that $(I,J)$ is an asymptotic pair for $(C(X),C(f))$, since $f$ preserve the orientation and $(I,J)$ is proximal for $(C(X),C(f))$. Now if $I=[a,b]$ is an arc proximal to the whole circle $X$. In this case $(a,b)$ is an asymptotic pair for $(X,f)$. Again as $f$ preserve the orientation and $(I,X)$ is proximal. It follow immediately that the sequence $(f^{n}(I))_{n\geq 0}$ converges to $X$ with respect to the Hausdrooff metric.

\textit{Case.2} $f$ reverse the orientation: Then $f^{2}$ preserve the orientation. Recall that a pair $(x,y)$ is a Li-Yorke pair with respect to a dynamical system $(Z,g)$ if and only if it is so for $(Z,g^{2})$. Thus by case $1$ above, $(C(X),C(f))$ has no Li-Yorke pair.
\end{proof}

\textbf{\textit{Proof of Theorem \ref{Main 1 eq}.}}
If $X$ is simple closed curve then property (iv) obviously holds (since $B\setminus A$ has at most two connected components for any subcontinua $A$ and $B$). We will now verify property (iii). Take any periodic subcontinuum $A$ of $X$ and let $C$ be a connected component of $X\setminus \Omega(f)$. Suppose that $\emptyset\neq A\cap C\varsubsetneq C$. In this case $A$ is an arc with periodic endpoints, at least one of them is in $C$, a contradiction. Thus $C\subset A$. Therefore property (iii) holds.   For properties (i) and (ii), see \cite{Lambard} and Lemma \ref{cercle}.

If $X$ is not a simple closed curve. The proof of $(i)\Longrightarrow (iv)$: Suppose that there exist periodic subcontinua $A$ and $B$ such that $A\subset B$ and some connected component of $B\setminus A$ has infinite orbit. We may assume that the $B$ and $A$ are fixed points of $(C(X),C(f))$ (otherwise consider a convenient iterate of $f$ instead). By collapsing $A$ to a single point $a$, we get a factor system $(\widetilde{X},\widetilde{f})$ of $(X,f)$ throughout the monotone natural quotient map $\pi:X\to \widetilde{X}$. Let $S^{\prime}=\{a\} \cup \big( \displaystyle \bigcup_{n\in \mathbb{Z}}\pi(f^{n}(I))\big)$, it is an $\widetilde{f}$-invariant subset of $\widetilde{X}$ and moreover $(S^{\prime},\widetilde{f}_{\mid S^{\prime}})$ is conjugated to the system $(S,g)$ defined in \cite{IH}. Hence $h_{top}(C(\widetilde{f}_{\mid S^{\prime}}))=h_{top}(C(g))$. By the facts that $h_{top}(C(g))=+\infty$ (see \cite{IH}) and $h_{top}(C(\widetilde{f}))\geq h_{top}(C(\widetilde{f}_{\mid S^{\prime}}))$, we get $h_{top}(C(\widetilde{f}))=+\infty$. It follows by Lemma \ref{factor factor} that $h_{top}(C(f))=+\infty$ (since the factor map $\pi:X\to \widetilde{X}$ is monotone).

Proofs of the implications $(iv)\Longrightarrow (iii)$, $(iii)\Longrightarrow (ii)$ and $(ii)\Longrightarrow (i)$ come respectively from Lemmas \ref{L1,premier equiv}, \ref{L44} and [Corollary  2.4 in \cite{LYP}]. The implication $(v)\Longrightarrow (ii)$ is clear and $(iii)\Longrightarrow (v)$ comes from Lemma \ref{L44}. \qed

\newpage

\begin{example}\label{example rational}
 \rm{\textit{The space $Y$:} We will define $Y$ as a subset of the real plane $\mathbb{R}^2$ as follow:
First let $I$ be the segment in $\mathbb{R}^2$ joining $(0,0)$ and $(0,1)$. Let $I_1(1)$  be the segment in $\mathbb{R}^2$ joining $(0,0)$ and $(1,1)$ and let $I_2(1)$ be the symmetric image of $I_1(1)$ with respect to the x-axis.

For each integer $n>1$ and each $i\in \{1,\dots,n\}$, let $I_{i}(n)$ be the segment in $\mathbb{R}^2$ joining $(0,0)$ and $\displaystyle{(\frac{i}{n},\frac{1}{n-i+1+\sum_{k=1}^{n-1}k})}$ and let $I_{n+i}(n)$  be the symmetric image of the segment $I_{n-i+1}(n)$ with respect to the x-axis. Set
$$Y=I\cup \big(\bigcup_{n\in\mathbb{N}}\cup_{1\leq i\leq 2n} I_{i}(n))\big).$$}

%
\end{example}

It is easy to check that any point in $I$ has a basis system of neighbourhoods with countable boundaries and any point outside $I$ is of order $\leq 2$. So $Y$ is a rational curve.

\textit{The map $F$:} First, let any point in $I$ be fixed by $F$. For each integer $n>1$ and each $i\in\{1,\dots,2n\}$, let $F_{\mid I_i(n)}$ be the similitude of the plane restricted to the segment $I_i(n)$ and sending $I_i(n)$ to $I_{i+1 \ mod(2n)}(n)$. In this way $F$ is obviously bijective. Observe that $F$ is locally a similitude on $Y\setminus I$, thus $F$ is continuous on $Y\setminus I$. Now observe that the continuity on $I$ comes from the fact that the speed of the orbit of $I_1(n)$ decreases to $0$ as $n\to\infty$. It turns out that $F$ is a self homeomorphism of $Y$.

\textit{Computation of topological entropy:}  Let $k,m\in\mathbb{N}$.  For each $\sigma=(\sigma_{1},\dots, \sigma_{2m}) \in \{1,\dots,k\}^{2m}$, we will define the subcontinuum $A_{\sigma}$ as follow: For each $i\in\{1,\dots,2m\}$, let $l_i(2m)$ denote the length of the segment $I_i(2m)$, then let $J_i^{\sigma}$  be the subarc of $I_i(2m)$ of length $\frac{\sigma_i}{k}l_i(2m)$ having $(0,0)$ as an endpoint. Thus define
$$A_{\sigma}=\displaystyle\cup_{1\leq i\leq 2m} J_i^{\sigma}.$$

It is easy to verify that $\{A_{\sigma}, \ \sigma\in \{1,\dots,k\}^{2m}\}$ is an $(2m,C(F),\frac{1}{k})$-separated set. It follows that
$sep(2m, C(F), \frac{1}{k})\geq k^{2m}$. Then $$\limsup_{m\to\infty}\frac{1}{m}\ln(sep(m,C(F),\frac{1}{k}))\geq \ln(k).$$ In result, $h_{top}(C(F))=h_{top}(2^F)=\infty$.

\bibliographystyle{unsrt}

\end{document}